\documentclass[12pt,a4paper]{article}
\usepackage{amsfonts}
\usepackage{bbm}

\usepackage{amssymb}

\usepackage{amsmath}

\newtheorem{theorem}{Theorem}[section]

\newtheorem{lemma}[theorem]{Lemma}

\newtheorem{example}[theorem]{Example}

\begin{document}
\title{Gr\"{o}bner-Shirshov bases for Coxeter
groups I\footnote{Supported by the NNSF of China (No.10771077) and
the NSF of Guangdong Province (No.06025062).} }
\author{
 Yuqun Chen  and Cihua Liu\\
{\small \ School of Mathematical Sciences, South China Normal
University}\\
{\small Guangzhou 510631, P. R. China}\\
{\small  yqchen@scnu.edu.cn}\\
{\small langhua01duo@yahoo.com.cn}}
\date{}
\maketitle \noindent\textbf{Abstract:} A conjecture of
Gr\"{o}bner-Shirshov basis of any Coxeter group has proposed  by
L.A. Bokut and L.-S. Shiao  \cite{bs01}.  In this paper, we give an
example to show that the conjecture is not true in general.  We list
all possible nontrivial inclusion compositions when we deal with the
general cases of the Coxeter groups. We give a Gr\"{o}bner-Shirshov
basis of a Coxeter group which is without nontrivial inclusion
compositions mentioned the above.

\noindent \textbf{Key words: } Coxeter group, Gr\"{o}bner-Shirshov
basis, word problem.

\noindent {\bf AMS} Mathematics Subject Classification(2000):
 20F05, 20F10, 20F55, 16S15, 13P10

\section{Introduction}

Let $M=\|m_{ij}\|_{n\times n}$ be a symmetric $n\times n$ matrix
such that $m_{ii}=1,\ 2\leq m_{ij}\leq\infty$. The Coxeter group
$W=W(M)$ is defined by the generators $s_1,\cdots, s_n$ and the
defining relations $(s_is_j)^{m_{ij}}=1$.

A conjecture of Gr\"{o}bner-Shirshov basis of any Coxeter group has
proposed by L.A. Bokut and L.-S. Shiao \cite{bs01}.
Gr\"{o}bner-Shirshov bases of all finite Coxeter groups were given
in \cite{bs01, Lee, Sv}.  As it is hypothesis, the conjecture is
true for any finite Coxeter group. In this paper, we give an example
to show that the above conjecture is not true in general.  We list
all possible nontrivial inclusion compositions (four cases) when we
deal with the general cases of the Coxeter groups. We then give a
new conjecture and prove it is true in some cases. We give a
Gr\"{o}bner-Shirshov basis of a Coxeter group which is without
nontrivial inclusion compositions mentioned the above. We give some
examples of such Coxeter groups  but not the finite Coxeter groups.
 We will consider other cases in
another papers in the future.

\section{Preliminaries}

We first cite some concepts and results from the literature
\cite{Sh, b72, b76} which are related to Gr\"{o}bner-Shirshov bases
for  associative algebras. A notion of the pre-Gr\"{o}bner-Shirshov
basis is new.

Let $X$ be a set and $F$ a field,  $F\langle X\rangle$ the free
associative algebra over $F$ generated by $X$, and $X^*$ the free
monoid generated by $X$. A well ordering $<$ on $X^*$ is monomial if
for any $u, v\in X^*$,
$$
u < v \Rightarrow w_{1}uw_{2} < w_{1}vw_{2},  \ for \  all \
 w_{1}, \ w_{2}\in  X^*.
$$

For any $u\in X^*$, denote by $|u|$ the length of $u$.

A standard example of monomial ordering on $X^*$ is the deg-lex
ordering which first compare two words by length and then by
comparing them lexicographically, where $X$ is a well ordered set.

Then, for any polynomial $f\in F\langle X\rangle$, $f$ has the
leading (maximal) word $\overline{f}$. We call $f$ {\it monic} if
the coefficient of $\overline{f}$ is 1.

Let $f,\ g\in F\langle X\rangle$ be two monic polynomials and $w\in
X^*$.

If  $w=\overline{f}b=a\overline{g}$  for some $a,b\in X^*$ such that
$|\overline{f}|+|\overline{g}|>|w|$, then $(f,g)_w=fb-ag$ is called
the {\it intersection composition }of $f,g$ relative to $w$.

If $w=\overline{f}=a\overline{g}b$  for some $a, b\in X^*$, then
$(f,g)_w=f-agb$ is called the {\it inclusion composition} of $f,g$
relative to $w$. The transformation $f\mapsto f-agb$ is called the
elimination of leading word (ELW) of $g$ in $f$.

In $(f,g)_w$, $w$ is called the {\it ambiguity} of the composition.

Let $S\subset F\langle X\rangle$ be a monic set. A composition
$(f,g)_w$ is called trivial modulo $(S,w)$, denoted by
$$
(f,g)_w\equiv0 \ \ \ mod(S,w)
$$
if $(f,g)_w=\sum\alpha_ia_is_ib_i,$ where every $\alpha_i\in F, \
s_i\in S,\ a_i,b_i\in X^*$, and $a_i\overline{s_i} b_i<w$.

Generally, for $f,g\in F\langle X\rangle,\ f\equiv g \ \ \ mod(S,w)$
we mean $f-g=\sum\alpha_ia_is_ib_i,$ where every $\alpha_i\in F, \
s_i\in S,\ a_i,b_i\in X^*$, and $a_i\overline{s_i} b_i<w$.

Recall that $S$ is a {\it Gr\"{o}bner-Shirshov basis} if any
composition of polynomials from $S$ is trivial modulo $S$.

\ \

Let $f$ and $r_1$ be two polynomials.  Then $f\mapsto f_1$ by ELW of
$r_1$ in $f$ means $f=\alpha_1a_1r_1b_1+f_1$ where $a_1,b_1\in X^*,\
\alpha_1\in F$ and $\bar f=a_1\overline{r_1} b_1$. Generally,
$f\mapsto f_1\mapsto\cdots \mapsto f_n\mapsto r$ means that $f=\sum
\alpha_ia_ir_ib_i+r$ where $\bar f=a_1\overline{r_1}
b_1>a_2\overline{r_2}b_2>\cdots
>a_n\overline{r_n}b_n>r$. If this is the case, we say that $f$ can be reduced to
$r$ via $\{r_1,\dots, r_n\}$.

Clearly, if $(f,g)_w$ can be reduced to zero by ELW of $S$, then
$(f,g)_w\equiv 0\ \ mod(S,w)$.

 The following lemma was first proved
by Shirshov \cite{Sh} for  free Lie algebras (with deg-lex ordering)
(see also Bokut \cite{b72}). Bokut \cite{b76} specialized the
approach of Shirshov to associative algebras (see also Bergman
\cite{b}). For commutative polynomials, this lemma is known as
Buchberger's Theorem (see \cite{bu65, bu70}).

\begin{lemma}\label{l1}
{\em (Composition-Diamond Lemma)} \ Let $F$ be a field, $A=F\langle
X|S\rangle=F\langle X\rangle/Id(S)$ and $<$ a monomial ordering on
$X^*$, where $Id(S)$ is the ideal of $F\langle X\rangle$ generated
by $S$. Then the following statements are equivalent:

\begin{enumerate}
\item[(1)] $S$ is a Gr\"{o}bner-Shirshov basis.
\item[(2)] $f\in Id(S)\Rightarrow \bar{f}=a\bar{s}b$
for some $s\in S$ and $a,b\in  X^*$.
\item[(\ref{e3})] $Irr(S) = \{ u \in X^* |  u \neq a\bar{s}b,s\in S,a ,b \in X^*\}$
is a $F$-basis of the algebra $A=F\langle X| S \rangle$.
\end{enumerate}

\end{lemma}

If a subset $S$ of $F\langle X\rangle$ is not a Gr\"{o}bner-Shirshov
basis then one can add all nontrivial compositions of polynomials of
$S$ to $S$. Continuing this process repeatedly, we finally obtain a
Gr\"{o}bner-Shirshov basis $S^{comp}$ that contains $S$. Such a
process is called Shirshov algorithm.

A set $S$ is called {\it reduced Gr\"{o}bner-Shirshov basis} if it
is a Gr\"{o}bner-Shirshov basis and there are no inclusion
compositions in $S$.

A set $S$ is called {\it pre-Gr\"{o}bner-Shirshov basis} if there
exists a subset $R\subset F\langle X\rangle$ such that the following
conditions hold.

(i)  $Id(R)=Id(S)$ and $R$ is a Gr\"{o}bner-Shirshov basis.  $R$ is
called a Gr\"{o}bner-Shirshov basis with related to $S$.

(ii)   For any $r\in R$, there exists $s\in S$ with $|\bar s|=|\bar
r|$ such that either $r=s$ or there exists a finite sequence of
ELW's of $S\setminus\{s\}$, $ s=s_0\mapsto s_1\mapsto\cdots \mapsto
s_n=r$, i.e., $s$ can be reduced to $r$ via $S\setminus\{s\}$.

\begin{lemma}
Let $S\subset F\langle X\rangle$ be an effective set (in a plurally
algebraic language, one may say that for any $n\geq 0$, one knows
all polynomials $s\in S_n$ of degree less or equal $n$ from $S$, and
there are finite number of these polynomials.) If $S$ is a
pre-Gr\"{o}bner-Shirshov basis, then the word problem is solvable
for the algebra $F\langle X| S\rangle=F\langle X\rangle/Id (S)$.
\end{lemma}
\textbf{Proof} Let $f\in F\langle X\rangle$ be a polynomial of
degree $n\geq 1$, $R$ be  a Gr\"{o}bner-Shirshov basis with related
to $S$. Then $f\in Id(S)$ iff $f$ goes to $0$ by the ELW of $R$. So
we need only to know all polynomials $r\in R_n$ of degree less or
equal than $n$ from $R$. From the definition of a
pre-Gr\"{o}bner-Shirshov basis, $R_n$ is a result of the ELW of
$S_n$ for polynomials from $S_n$. Since we know $S_n$, we can find
$R_n$ effectively. \hfill $\blacksquare$

 \ \

Let $A=sgp\langle X|S\rangle$ be a semigroup presentation. Then $S$
is also a subset of $F\langle X \rangle$ and we can find
Gr\"{o}bner-Shirshov basis
 $S^{comp}$. We also call $S^{comp}$ a
Gr\"{o}bner-Shirshov basis of $A$. The  set $Irr(S^{comp})=\{u\in
S^*|u\neq a\overline{f}b,\ a ,b \in X^*,\ f\in S^{comp}\}$ is a
linear basis of $F\langle X|S\rangle$ which is also a set of all
normal forms of $A$.

\section{Gr\"{o}bner-Shirshov bases of Coxeter groups}

Let $\Sigma=\{\sigma_1, \cdots, \sigma_n\}$ be a finite set. Let
$M=(m_{ij})$ be a symmetric $n\times n$ matrix over the natural
numbers together with $\infty$, such that $m_{ii}=1,\ 2\leq
m_{ij}\leq\infty$ for $i\neq j$. Such an M is called a Coxeter
matrix. Now, we use $W$ to denote
$$
W=W(M)=sgp\langle \Sigma|(\sigma_i\sigma_j)^{m_{ij}}=1,\ 1\leq i,\
j\leq n,\ m_{ij}\neq\infty\rangle.
$$
W is called the Coxeter group (see, for example, \cite{S02}) with
respect to Coxeter matrix $M$.

We order $\Sigma^*$ by the  deg-lex ordering, where
$\sigma_1<\dots<\sigma_n$.

For any $i,j\ \ (1\leq i,j\leq n)$, denote by
$m_{\sigma_{_i}\sigma_{_j}}=m_{ij}$. For any $s,s'\in \Sigma$, we
now define for finite $m_{ss'}$ the following notation:

\ \ \ \ \ \ $m(s,s')=ss'\cdots $(there are $m_{ss'}$ alternative
letters $s,s')$,

\ \ \ \ \ $(m-i)(s,s')=ss'\cdots$ (there are $m_{ss'}-i$ alternative
letters $s,s'$, $1\leq i\leq m_{ss'}$). \\
With the above notation, the  defining relations of W can be
presented in the following forms
\begin{eqnarray}
\label{e1}&&s^2=1 \\
\label{e2}&&m(s,s')=m(s',s),\ s>s'
\end{eqnarray}
for all $s,s'\in \Sigma$ and finite $m_{ss'}$.

Define $s\rhd s'$ if  $s>s'$ and $m_{ss'}=2$.

\begin{lemma} (\cite{bs01})\label{s.1}
In group $W$, we have
\begin{eqnarray}
\nonumber &&(m-1)(s_{0},s'_0)(m-1)(s_1,s'_1)\cdots
(m-1)(s_{k},s'_{k})
m(s_{k+1},s'_{k+1})\\
\label{e3}&=&m(s'_0,s_0)(m-1)(s_1,s'_1)\cdots
(m-1)(s_{k},s'_{k})(m-1)(s_{k+1},s'_{k+1})
\end{eqnarray}
where $k\geq 0,\ s_0,s'_0,\dots,s_{k+1}, s'_{k+1}\in \Sigma$
and for any $i, \ 0\leq i\leq k$

$s'_{i+1}= \left\{
\begin{array}{ll} s'_{i} \ \ \ \ \ \
 \mbox{ if }\ m_{s_is'_{i}}\ \mbox{ is
even},\\
s_{i} \ \ \ \ \ \   \mbox{ if }\ m_{s_{i}s'_{i}}\ \mbox{ is odd}.
\end{array}\right.$
\end{lemma}
\textbf{Proof}\ Since
\begin{eqnarray*}
&&(m-1)(s_{0},s'_0)(m-1)(s_1,s'_1)\cdots (m-1)(s_{k},s'_{k})
m(s_{k+1},s'_{k+1})\\
&=&(m-1)(s_{0},s'_0)(m-1)(s_1,s'_1)\cdots (m-1)(s_{k},s'_{k})
m(s'_{k+1},s_{k+1})\\
&=&(m-1)(s_{0},s'_0)(m-1)(s_1,s'_1)\cdots m(s_{k},s'_{k})
(m-1)(s_{k+1},s'_{k+1})\\
&=&\cdots\\
&=&m(s'_0,s_0)(m-1)(s_1,s'_1)\cdots
(m-1)(s_{k},s'_{k})(m-1)(s_{k+1},s'_{k+1}),
\end{eqnarray*}

we obtain the result. \hfill $\blacksquare$

\ \

Denote by
$$
S=\{(\ref{e1}),(\ref{e2}),(3')\}
$$
where $(3')$ consists of all relations in (\ref{e3}) with the extra
properties
\begin{eqnarray}
\label{e4}&&s_0>s'_0,\ s_1<s'_1,\ \cdots,\ s_k<s'_k,\ s_{k+1}<s'_{k+1}\\
\label{e5}&&\{s_i,s'_i\}\neq\{s_{i+1},s'_{i+1}\},\ 0\leq i\leq k
\end{eqnarray}

It was conjectured in \cite{bs01} that a Gr\"{o}bner-Shirshov basis
of $W$ can be obtained from $S$ using only commutative relations of
$W$ $(m(s,s')=m(s',s)$ where $m_{ss'}=2$). The following example
shows that this conjecture is not true in general.

\begin{example}\label{ex1}
Let $\Sigma=\{s_1,s_2,s_3,s_4\}$ with $s_1<s_2<s_3<s_4$,
$M=(m_{ij})$ the $4\times4$ Coxeter matrix where
$m_{s_1s_2}=m_{s_2s_3}=m_{s_2s_4}=\infty,\ m_{s_1s_3}=3,\
m_{s_1s_4}=2,\ m_{s_3s_4}=5$ and $m_{s_is_i}=1,\ i=1,2,3,4$. Then
\begin{eqnarray*}
(\ref{e1})&=&\{s_i^2=1,\ i=1,2,3,4\},\\
(\ref{e2})&=&\{s_4s_1=s_1s_4,\
s_3s_1s_3=s_1s_3s_1,\ s_4s_3s_4s_3s_4=s_3s_4s_3s_4s_3\},\\
(3')&=&\{(m-1)(s_4,s_3)m(s_1,s_4)=m(s_3,s_4)(m-1)(s_1,s_4),\\
&&\
(m-1)(s_4,s_3)(m-1)(s_1,s_4)m(s_3,s_4)=m(s_3,s_4)(m-1)(s_1,s_4)(m-1)(s_3,s_4),\\
&&\
 (m-1)(s_4,s_3)(m-1)(s_1,s_4)(m-1)(s_3,s_4)m(s_1,s_3)\\
 &&=m(s_3,s_4)(m-1)(s_1,s_4)(m-1)(s_3,s_4)(m-1)(s_1,s_3)\}.
\end{eqnarray*}
A Gr\"{o}bner-Shirshov basis of $W$ is $(\ref{e1})\cup
(\ref{e2})\cup (3'')$, where
\begin{eqnarray*}
(3'')&=&\{(m-1)(s_4,s_3)m(s_1,s_4)=m(s_3,s_4)(m-1)(s_1,s_4),\\
&&
(m-3)(s_4,s_3)s_1s_4s_3s_1(m-1)(s_4,s_3)=m(s_3,s_4)(m-1)(s_1,s_4)(m-1)(s_3,s_4),\\
&&
 (m-3)(s_4,s_3)s_1s_4s_3s_1(m-3)(s_4,s_3)s_1s_4(m-1)(s_3,s_1)\\
 &&=m(s_3,s_4)(m-1)(s_1,s_4)(m-1)(s_3,s_4)(m-1)(s_1,s_3)\}
\end{eqnarray*}
which are obtained from $(3')$ by using the relations
$s_4s_1=s_1s_4,\ s_3s_1s_3=s_1s_3s_1$. \hfill $\blacksquare$
\end{example}

\ \

Then we give the following conjecture.

\ \

\noindent\textbf{Conjecture (L.A. Bokut):} The set of relations
(\ref{e1}),(\ref{e2}),(\ref{e3}) is a pre-Gr\"{o}bner-Shirshov basis
of $W$.

\ \

In this paper, we will show that the above new conjecture is true
when $M$ satisfies some conditions.

\begin{theorem} \label{s.2}
Let $S=\{(\ref{e1}),(\ref{e2}),(3')\}$. Then if $S$ is a
pre-Gr\"{o}bner-Shirshov basis of $W$ then so is
$\{(\ref{e1}),(\ref{e2}),(\ref{e3})\}$.
\end{theorem}
\textbf{Proof} It suffices  to show that for any
\begin{eqnarray*}
f&=&(m-1)(s_{0},s'_0)(m-1)(s_1,s'_1)\cdots (m-1)(s_{k},s'_{k})
m(s_{k+1},s'_{k+1})\\
&&-m(s'_0,s_0)(m-1)(s_1,s'_1)\cdots
(m-1)(s_{k},s'_{k})(m-1)(s_{k+1},s'_{k+1})
\end{eqnarray*}
 in $(\ref{e3})$
without property $(\ref{e4})$ or $(\ref{e5})$, $f$ has an
expression: $f=\sum a_ir_ib_i$, where $r_i\in S,\ a_i,b_i\in X^*$.
We prove this by induction on $k$.

For $k=0$,
$f=(m-1)(s_0,s'_0)m(s_1,s'_1)-m(s'_0,s_0)(m-1)(s_1,s'_1)$. There are
two cases to consider.

Case 1. $f$ is without property $(\ref{e4})$.

If $s_1>s'_1$, then
$$
f=(m-1)(s_0,s'_0)(m(s_1,s'_1)-m(s'_1,s_1))+(m(s_0,s'_0)-m(s'_0,s_0))(m-1)(s_1,s'_1).
$$

If $s_0<s'_0$, then
$$
f=-(m(s'_0,s_0)-m(s_0,s'_0))(m-1)(s_1,s'_1)+(m-1)(s_0,s'_0)(m(s_1,s'_1)-m(s'_1,s_1)).
$$

Case 2. $f$ is without property $(\ref{e5})$.

If $\{s_0,s'_0\}=\{s_1,s'_1\}$, then by ELW's of $s_0^2=1$ and
$s_0'^2=1$, $f\mapsto \cdots \mapsto
f_{m_{_{s_0s_0'}-1}}=s'_0-s'_0=0$.

Thus the result is true for $k=0$.

For $k>0$, there are also two cases to consider.

Case 1. $f$ is without property $(\ref{e4})$.

If $s_{k+1}>s'_{k+1}$, then
$$
f=(m-1)(s_{0},s'_0)(m-1)(s_1,s'_1)\cdots
(m-1)(s_{k},s'_{k})r_1+r_2(m-1)(s_{k+1},s'_{k+1})
$$
where $r_1=m(s_{k+1},s'_{k+1})-m(s'_{k+1},s_{k+1})\in (\ref{e2})$
and
$$r_2=(m-1)(s_{0},s'_0)(m-1)(s_1,s'_1)\cdots
m(s_{k},s'_{k})-m(s'_0,s_0)(m-1)(s_1,s'_1)\cdots
(m-1)(s_{k},s'_{k})\in (\ref{e3}).$$ By induction, $r_2$ is a
combination of relations in $(3')$. Then the result follows.

If $s_0'>s_0$, then
$$
f=-r_1(m-1)(s_1,s'_1)\cdots (m-1)(s_{k},s'_{k})
(m-1)(s'_{k+1},s_{k+1})+(m-1)(s_0,s'_0)r_2
$$
where $r_1=m(s'_0,s_0)-m(s_0,s'_0)\in (\ref{e2})$ and
$
r_2=(m-1)(s_1,s'_1)\cdots (m-1)(s_{k},s'_{k})
m(s'_{k+1},s_{k+1})-m(s'_1,s_1)\cdots
(m-1)(s_{k},s'_{k})(m-1)(s_{k+1},s'_{k+1})
$
is in (\ref{e3}). By induction, the result follows.

If there exists $i,\ 0<i<k+1$ such that  $s_i>s'_i,\ s_0>s'_0,\
s_{k+1}<s'_{k+1}$, then
$$
f=(m-1)(s_{0},s'_0)\cdots
(m-1)(s_{i-1},s'_{i-1})r_1+r_2(m-1)(s_{i},s'_{i})\cdots
(m-1)(s_{k+1},s'_{k+1})
$$
where $r_1=(m-1)(s_i,s'_i)\cdots
m(s_{k+1},s'_{k+1})-m(s'_i,s_i)\cdots (m-1)(s_{k+1},s'_{k+1})$,
$r_2=(m-1)(s_{0},s'_0)\cdots m(s_{i-1},s'_{i-1})-m(s'_0,s_0)\cdots
(m-1)(s_{i-1},s'_{i-1})$, and both of them are in (\ref{e3}). By
induction, the result follows.

Case 2. $f$ is without property $(\ref{e5})$.

Let us have $f$ with condition (\ref{e4}). Suppose
$\{s_i,s'_i\}=\{s_{i+1},s'_{i+1}\},\ 0\leq i\leq k.$

If $i<k$, then by ELW's of  $s_i^2=1$ and $s_i'^2=1$,
\begin{eqnarray*} f\mapsto\cdots &\mapsto&(m-1)(s_0,s'_0)\cdots
(m-1)(s_{i-1},s'_{i-1})(m-1)(s_{i+2},s'_{i+2})\cdots
m(s_{k+1},s'_{k+1})\\
&&-m(s_0',s_0)\cdots
(m-1)(s_{i-1},s'_{i-1})(m-1)(s_{i+2},s'_{i+2})\cdots(m-1)(s_{k+1},s'_{k+1})
\end{eqnarray*}
is in (\ref{e3}) since $s'_{i+2}$ is the last second letter of
$(m-1)(s_{i+1},s'_{i+1})$ which, in fact, is $s'_i$. By induction,
the result follows.

If $i=k$ then by ELW's of  $s_k^2=1$ and $s_k'^2=1$,
\begin{eqnarray*}
f\mapsto\cdots &\mapsto&(m-1)(s_0,s'_0)\cdots
m(s_{k-1},s'_{k-1})-m(s'_0,s_0)\cdots (m-1)(s_{k-1},s'_{k-1})
\end{eqnarray*}
is in (\ref{e3}). By induction, the result follows. \hfill
$\blacksquare$

\ \

We will deal with inclusion compositions $(f,g)_w,\ \bar f=a\bar gb,
\ w=\bar f$ and $f\in (3'), \ g\in (2)\cup(3')$. We will prove that
in the most cases they are trivial except six cases in Theorems
\ref{s.11},  \ref{s.14}, \ref{s.16} and  \ref{s.17}.

\ \

\noindent {\bf Notation}:

\ \

We will fix two ``typical" relations in $(3')$.

Let $f$ be a  relation in $(3')$,
\begin{eqnarray}\label{e6}
&&f=u_0u_1\cdots u_ku_{k+1}y_{k+1}-s_0'u_0u_1\cdots u_ku_{k+1}=\bar f-f_0\\
\nonumber && u_i=(m-1)(s_i,s'_i),\\
\nonumber &&x_{i}\ \mbox{ the\  last\  letter\  of}\  (m-1)(s_i,s'_i),\\
\nonumber &&y_{i}\ \mbox{ the\  last\  letter\  of}\   m(s_i,s'_i),\
\ \ \ \ \ \ \ \ \ \  0\leq i\leq k+1
\end{eqnarray}
where
$
\{x_i,s'_{i+1}\}=\{s_i,s'_i\},\ y_i=s'_{i+1},\
m(s_is'_i)=(m-1)(s_i,s'_i)s'_{i+1},\ \ \ \ \ 0\leq i\leq k.
$

Let $g$ be an other relation in $(3')$,
\begin{eqnarray}\label{e7}
&&g=v_0v_1\cdots v_qv_{q+1}z_{q+1}-p'_0v_0v_1\cdots v_qv_{q+1}=\bar g-g_0\\
\nonumber && v_i=(m-1)(p_i,p'_i),\\
\nonumber &&t_{i}\ \mbox{ the\  last\  letter\  of}\  v_i,\\
\nonumber &&z_{i}\ \mbox{ the\  last\  letter\  of}\   m(p_i,p'_i),\
\ \ \ \ \ \ \ \ \ \  0\leq i\leq q+1
\end{eqnarray}
where $\{t_i,p'_{i+1}\}=\{p_i,p'_i\},\ z_i=p'_{i+1},\
m(p_ip'_i)=(m-1)(p_i,p'_i)p'_{i+1},\ \ \ \ \  \ 0\leq i\leq q$.

\ \

In Lemmas (Theorems) \ref{s.3}--\ref{s.14}, we always assume that
$f,g\in (3')$ with the forms  (\ref{e6}),  (\ref{e7}) respectively
and $\bar f=a\bar gb$ for some words $a,b$.

\begin{lemma}\label{s.3}
If $\bar f=a\bar g$, then $a=1$ and $f=g$.
\end{lemma}
\textbf{Proof}\ Since $y_{k+1}=z_{q+1}$ and $x_{k+1}=t_{q+1}$,
$u_{k+1}y_{k+1}=v_{q+1}z_{q+1}$. Since $x_{k}=t_{q}$ and
$y_k=s'_{k+1}=p'_{q+1}=z_{q}$, $u_{k}y_{k}=v_{q}z_{q}$. Similarly,
we have $u_{k-1}y_{k-1}=v_{q-1}z_{q-1},\cdots,u_0y_0=v_0z_0$. Then
$a=1$ and $\bar f=\bar g$.

Noting that $u_0\cdots u_{k+1}=v_0\cdots v_{q+1}$, in order to prove
$\bar f=\bar g$ it is sufficient to show that $s'_0=p'_0$. Induction
on $k$.

If $k=0$, then $y_{1}=z_{q+1}$ and $x_{1}=t_{q+1}$. Then
$u_{1}y_{1}=v_{q+1}z_{q+1}$. Since $x_0=t_q$ and $s'_1=p'_{q+1}$,
$u_{0}y_0=v_qp'_{q+1}$. Then $q=k=0$ and $s'_0=p'_0$.

For $k>0$, we have $y_{k+1}=z_{q+1}$ and $x_{k+1}=t_{q+1}$,
$u_{k+1}y_{k+1}=v_{q+1}z_{q+1}$. Then $y_k=z_q$.

Let $h=u_0\cdots u_{k}y_k-s'_0u_0\cdots u_{k}$ and $ q= v_0\cdots
v_{q}z_{q}-p'_0v_0\cdots v_{q}$. Clearly, $\bar h=\bar q$, Then by
induction, we have $s'_0=p'_0$. \hfill $\blacksquare$

\begin{lemma}\label{s.4}
If there exist $i,j$ such that $s_i=p_j$, $s_i'=p_j'$ and $u_i$ is a
subword of $\bar g$, then $\bar f=\bar g$.
\end{lemma}
\textbf{Proof}\ If $i=0$ then $j=0$ since $u_iy_i=v_jz_{j}$. Then
$s'_1=y_0=z_0=p'_1$. Since $\bar g$ is a subword of $\bar f$,
$s_1=p_1$ and $u_2y_2=v_2z_2$. Hence $u_iy_i=v_iz_{i}$ for any $i,\
1\leq i\leq k+1$. Then $\bar f=\bar g$.

If $i\neq 0$, then $j\neq 0$. Otherwise, we have
$p_0=s_i<s'_{i+1}=p'_0$, a contradiction. Then $x_{i-1}=t_{j-1}$.
Since $y_i=s'_{i+1}=p'_{j+1}=z_j$, $u_{i-1}y_i=v_{j-1}z_{j}$.
Similarly, we have
$u_{i-2}y_{i-2}=v_{j-2}z_{j-2},\cdots,u_0y_0=v_0z_0$ and $j=i$.
Also, $s_{i+1}=p_{i+1},\ s'_{i+1}=y_{i}=z_{i}=p'_{i+1}$ imply that
$u_{i+1}y_{i+1}=v_{i+1}z_{i+1}$. Therefore,
$u_{i+2}y_{i+2}=v_{i+2}z_{i+2},\cdots,u_{k+1}y_{k+1}=v_{k+1}z_{k+1}$.
\hfill $\blacksquare$

\ \

In what follows we assume that $\bar f\neq \bar g$.

\begin{lemma}\label{ls.1}
If there exists $i>0$ such that $\ |u_i|>1$, $\bar g=cu_id$, $\bar
f=acu_idb$, $ac=u_0\cdots u_{i-1}$ and $c=v_0\cdots v_{j-1}$, then
$u_i=v_jv_{j+1}\cdots v_n$ and $|v_j|=\cdots=|v_n|=1$.

Moreover, if $u_{i+1}$ is also a subword of $\bar g$,  then $
u_{i+1}=v_{n+1}\cdots v_{l}$ such that $|v_j|=\cdots=|v_{l}|=1$.
\end{lemma}
\textbf{Proof} By Lemma \ref{s.4} and $\bar f\neq \bar g$, we have
$u_i\neq v_j$.

Since $v_j=(m-1)(s_i,p'_{j})$ and $u_i\neq v_j$, $p'_j\neq s'_{i}$.
Then $|v_j|=1$ and $v_{j+1}=(m-1)(s'_{i},p'_j)$. If $|v_{j+1}|>1$,
then $p'_j=s_{i+1}$ and $u_i=s_is'_{i}$. Now,
$s'_{i}<p'_j=s_{i+1}<s'_{i+1}=s_i$,  a contradiction. Then
$|v_{j+1}|=1$. This shows that $u_i=v_jv_{j+1}\cdots v_n$ such that
$|v_j|=\cdots=|v_n|=1$.

If $u_{i+1}$ is also a subword of $\bar g$, we have
$v_{n+1}=(m-1)(s_{i+1},p'_j)$. If $|u_{i+1}|>1$, then by a similar
proof of the above, we have $u_{i+1}=v_{n+1}\cdots v_{l}$ such that
$|v_{n+1}|=\cdots =|v_{l}|=1$. If $|u_{i+1}|=1$ and $|v_{n+1}|>1$,
then $p'_j=s_{i+2}$,
$s'_{i+1}<s_{i+2}<s'_{i+2}\in\{s_{i+1},s'_{i+1}\}$, a contradiction.
Therefore, $|v_{n+1}|=1$ and $u_{i+1}=v_{n+1}$. \hfill
$\blacksquare$

\ \

\begin{lemma}\label{ls.2}
If there exist $i,i'\ (i'\geq 1)$ such that $u_i\cdots
u_l=v_{i'}\cdots v_{q+1}$, then $|u_i|=\cdots =|u_l|=1$.
\end{lemma}
\textbf{Proof}\ Suppose there exists a minimal $ j\ (i\leq j\leq l)$
such that $|u_j|>1$. We will show that $\bar g=cu_jd$, where
$c=v_0\cdots v_n,\ i'-1\leq n\leq q$. Otherwise, $s_j$ is a subword
of $v_n$. Then $v_n=(m-1)(s_{j-1},s_j)=s_{j-1}s_j$ and
$v_{n+1}=(m-1)(s'_{j},s_{j-1})\ (j> 1)$. So, $s'_{j}<s_{j-1}<s_j$, a
contradiction.

Then by Lemma \ref{ls.1}, we have $u_j=v_{n+1}\cdots v_{l'}$ such
that $|v_{n+1}|=\cdots=|v_{l'}|=1$.

Moreover, $u_{j+1}\cdots u_l=v_{l'+1}\cdots v_{q+1}$ such that
$|v_{l'+1}|=\cdots=|v_{q+1}|=1$. Then $z_{q+1}=s_{l+1}$ and there
exists $v_p\ (n+1\leq p\leq q+1)$ such that $s'_{l+1}=v_p<s_{l+1}$,
a contradiction. \hfill $\blacksquare$

 \ \

\begin{lemma} \label{s.7}
If $\bar f=\bar gb$ with $b\neq 1$, then $|u_0|=1$ or $|u_0|=2$.
\end{lemma}
\textbf{Proof} If $|u_0|>2$, then $|v_0|=1$. Otherwise, by Lemma
\ref{s.4}, $\bar f=\bar g$, a contradiction. Clearly, $|v_1|=1$.
Then $p_2=s_{0}$ and $p_2<p'_2=p'_0<p_0=s_{0}$, a contradiction.
\hfill $\blacksquare$

\begin{lemma} \label{s.8}
Suppose that $\bar f=\bar gb=\bar
g(m-2)(s'_{l+1},s_{l+1})u_{l+2}\cdots u_{k+1}y_{k+1}$. Then
$|u_1|=\cdots =|u_l|=1$, $|v_0|=1$ and $(f,g)_{\bar f}\equiv 0$.
\end{lemma}
\textbf{Proof} There are two cases to consider.

Case 1. $u_0=s_0$.

We will show that $v_0=s_{0}$. Otherwise, $v_0=s_{0}s_1$. If
$|u_1|>1$, then $v_1=(m-1)(s'_1,s_0)=(m-1)(s'_0,s_0)=s'_0=s'_1$,
$u_1=s_1s'_1$ and $v_2\cdots v_{q+1}=u_2\cdots u_l$. By Lemma
\ref{ls.2}, we have $|u_2|=\cdots=|u_l|=1$. Then there exists
$s_{j}\in\{s_2,\cdots,s_{l+1}\}$ such that $s_{j}=s_{0}$, a
contradiction. Then $|u_1|=1$ and $v_1\cdots v_{q+1}=u_2\cdots u_l$.
By Lemma \ref{ls.2}, we have $|u_2|=\cdots=|u_l|=1$. This implies
that there exists $l+1\geq j>1$ such that $s_{j}=s_0$, a
contradiction.

Since $v_0=s_0$ and $v_1\cdots v_{q+1}=u_1\cdots u_l$, by Lemma
\ref{ls.2}, we have $|u_1|=\cdots=|u_l|=1$.

Suppose $p'_0=s_{j}$ where $s_{j}\in\{s_2,\cdots,s_{l+1}\}$. If
$j<l+1$, there exists an $i$ such that $|v_i|>1$ and so
$|u_{l+1}|=1$. By Lemma \ref{s.6}, $(f,g)_{\bar f}\equiv 0$.

If $j=l+1$, then $s_0\rhd s_{l+1}\rhd s_j$, $s_1\rhd s'_{l+1}\rhd
s_j$ for any $j,\ 1\leq j\leq l$ and
\begin{eqnarray*}
(f,g)_{\bar f}&\equiv& s_{l+1}s_0s_1\cdots
s_l(m-2)(s'_{l+1},s_{l+1})\cdots
m(s_{k+1},s'_{k+1})\\
&&-s'_0s_{l+1}s_0s_1\cdots s_l(m-2)(s'_{l+1},s_{l+1})\cdots
(m-1)(s_{k+1},s'_{k+1})\\
&\equiv&s_{l+1}s'_{l+1}s_0s_1\cdots
s_l(m-3)(s_{l+1},s'_{l+1})\cdots m(s_{k+1},s'_{k+1})\\
&&-s'_0s_{l+1}s'_{l+1}s_0s_1\cdots s_l(m-3)(s_{l+1},s'_{l+1})\cdots
(m-1)(s_{k+1},s'_{k+1})\\
&&\cdots \\
&\equiv& (m-1)(s_{l+1},s'_{l+1})s_0s_1\cdots
s_l(m-1)(s_{l+2},s'_{l+2})\cdots
 m(s_{k+1},s'_{k+1})\\
&&-s'_0(m-1)(s_{l+1},s'_{l+1})s_0s_1\cdots
s_l(m-1)(s_{l+2},s'_{l+2})\cdots(m-1)(s_{k+1},s'_{k+1})\\
&\equiv& (m-1)(s_{l+1},s'_{l+1})s'_{l+2}s_0s_1\cdots
s_l(m-1)(s_{l+2},s'_{l+2})\cdots
(m-1)(s_{k+1},s'_{k+1})\\
&&-(m-1)(s_{l+1},s'_{l+1})s'_{l+2}s_0s_1\cdots
s_l(m-1)(s_{l+2},s'_{l+2})\cdots (m-1)(s_{k+1},s'_{k+1})\\
&\equiv&0
\end{eqnarray*}
since $s'_{l+1}=\cdots =s'_0$,
$(m-1)(s_{l+1},s'_{l+1})s'_{l+2}=m(s_{l+1},s'_{l+1})$ and
$h=s_0s_1\cdots s_lu_{l+2}\cdots
 u_{k+1}y_{k+1}-s'_{l+2}s_0s_1\cdots s_lu_{l+2}\cdots
 u_{k+1}$ in (\ref{e3}) with property (\ref{e4}).

Case 2. $u_0=s_0s'_0$. Then $v_0=s_0,\ v_1=(m-1)(s'_0,p'_0)$. There
are two subcases to consider.

Subcase 1. $|v_1|>1$. Then $p'_0=s_1$ and $|u_1|=1$. If $|v_1|>2$,
then $s_2=s'_0,\ v_1=s'_0s_1s'_0$ and $u_2=s_2s'_2=s'_0s_0$. This
shows $v_2=(m-1)(s_0,s_1)$, a contradiction. Then $v_1=s'_0s_1$ and
$v_2\cdots v_{q+1}=u_2\cdots u_l$. By Lemma \ref{ls.2},
$|u_2|=\cdots=|u_l|=1$. Clearly, $s'_0\not\in
\{s_{2},\cdots,s_{l-1}\}$, otherwise,  there exists $u_i\ (2\leq
i\leq l)$ such that $s_{i-1}=s'_0$ and $u_{i}=(m-1)(s'_0,s_0)$ which
contradicts  $|u_i|=1$.

Then $|v_2|=\cdots=|v_{q+1}|=1$ and $s_{l+1}=s'_0$, $s'_{l+1}=s_0$.
Now,
\begin{eqnarray*}
(f,g)_{\bar f}&\equiv &s_{1}s_0s'_0s_1\cdots
s_ls'_{l+1}u_{l+2}\cdots u_{k+1}y_{k+1}-s'_0s_{1}s_0s'_0s_1\cdots
s_ls'_{l+1}u_{l+2}\cdots
u_{k+1}\\
&\equiv&s_{1}s'_0s_0s'_0s_1\cdots s_lu_{l+2}\cdots
u_{k+1}y_{k+1}-s'_0s_{1}s'_0s_0s'_0s_1\cdots s_lu_{l+2}\cdots
u_{k+1}\\
&\equiv&s_{1}s'_0s_0s_1s'_0s_1\cdots s_lu_{l+2}\cdots
u_{k+1}-s'_0s_{1}s'_0s_0s'_0s_1\cdots s_lu_{l+2}\cdots u_{k+1}
\end{eqnarray*}
since $h=s'_0s_1\cdots s_lu_{l+2}\cdots
u_{k+1}y_{k+1}-s_1s'_0s_1\cdots s_lu_{l+2}\cdots u_{k+1}$ in
(\ref{e3}) with property (\ref{e4}) and $s'_{l+2}=s_{l+1}=s'_0$.

Since $s_0\rhd s_1$, $m_{s_1s'_0}=3$ and $s_1>s'_0$, we have
$s_{1}s'_0s_0s_1\mapsto s_{1}s'_0s_1s_0\mapsto s'_0s_1s'_0s_0$ and
hence $(f,g)_{\bar f}\equiv 0$.

Subcase 2. $|v_1|=1$. Then $v_2\cdots v_{q+1}=u_1\cdots u_l$. By
Lemma \ref{ls.2}, we have $|u_1|=\cdots=|u_l|=1$.

Suppose $p'_0=s_{j}$ where $s_{j}\in \{s_1,\cdots,s_{l+1}\}$. If
$j<l+1$, then $|u_{l+1}|=1$. If $j=l+1$, we have $|u_{l+1}|=1$ since
$s_0\rhd s_{l+1}$. Then
\begin{eqnarray*}
(f,g)_{\bar f}&\equiv &s_{j}s_0s'_0s_1\cdots s_lu_{l+2}\cdots
u_{k+1}y_{k+1}-s'_0s_{j}s_0s'_0s_1\cdots s_lu_{l+2}\cdots
u_{k+1}\\
&\equiv&s_{j}s'_0s_0s'_0s_1\cdots s_lu_{l+2}\cdots
u_{k+1}-s'_0s_{j}s_0s'_0s_1\cdots s_lu_{l+2}\cdots u_{k+1} \\
&\equiv&s'_0s_{j}s_0s'_0s_1\cdots s_lu_{l+2}\cdots
u_{k+1}-s'_0s_{j}s_0s'_0s_1\cdots s_lu_{l+2}\cdots u_{k+1}\\
&\equiv& 0
\end{eqnarray*}
since $s_0s'_0s_1\cdots s_lu_{l+2}\cdots
u_{k+1}y_{k+1}-s'_0s_0s'_0s_1\cdots s_lu_{l+2}\cdots u_{k+1}$ is in
(\ref{e3}). \hfill $\blacksquare$

\ \

\begin{lemma} \label{s.6}
If $|u_i|=\cdots =|u_{l+1}|=1$ and $u_i\cdots u_{l+1}=\bar g$, then
$(f,g)_{\bar f}\equiv 0$.
\end{lemma}
\textbf{Proof} Clearly,  $\bar g=u_i\cdots u_{l+1}\mapsto
u_ju_i\cdots u_{l}=g_0$ for some $i<j\leq{l+1}$.

If $i=0$, then
$$
(f,g)_{\bar f}\equiv u_ju_0\cdots u_lu_{l+2}\cdots
u_{k+1}y_{k+1}-u_js'_0u_i\cdots u_lu_{l+2}\cdots u_{k+1}\equiv s_jh
$$
where $h=u_0\cdots u_lu_{l+2}\cdots u_{k+1}y_{k+1}-s'_0u_0\cdots
u_lu_{l+2}\cdots u_{k+1}$ is  in (\ref{e3}) with property (\ref{e4})
and $s_j\bar h<\bar f$.  By Theorem \ref{s.2}, the result follows.

If $i>0$, then
$$
(f,g)_{\bar f}\equiv u_0\cdots u_{i-1}u_ju_i\cdots u_lu_{l+2}\cdots
u_{k+1}y_{k+1}-s_0u_0\cdots u_{i-1}u_ju_i\cdots u_lu_{l+2}\cdots
u_{k+1}\triangleq h
$$
where $h$ is  in (\ref{e3}) with property (\ref{e4}) and $\bar
h<\bar f$.  By Theorem \ref{s.2}, the result follows. \hfill
$\blacksquare$

 \ \

 The following lemmas are dealing with the case $\bar f=a\bar g b$,
 $a\neq 1,\ b\neq 1$.

 \ \

In Lemmas (Theorems) \ref{s.9}--\ref{s.15}, $i$ and $l$ are fixed
such that $0\leq i<l\leq k$, $u_0\cdots u_{i-1}=1$ if $i=0$ and
$u_{l+2}\cdots u_{k+1}=1$ if $l=k$.

\ \
 \begin{lemma} \label{s.9}
If $\bar f=u_0\cdots u_{i-1}(m-2)(s_i,s'_i)\bar
g(m-2)(s'_{l+1},s_{l+1})u_{l+2}\cdots u_{k+1} y_{k+1}$, then
$|u_{i+1}|=\cdots =|u_l|=1$.
\end{lemma}
\textbf{Proof} There are three cases to consider.

Case 1.  $v_0=x_{i}$. Then $v_1\cdots v_{q+1}=u_{i+1}\cdots u_{l}$.
By Lemma \ref{ls.2}, we have $|u_{i+1}|=\cdots=|u_l|=1$.

Case 2. $v_0=(m-1)(x_{i},s_{i+1})$ and $|v_0|>2$. Then $|u_{i+1}|=1$
and $s_{i+2}=x_{i}$. If $|u_i|>1$, then $v_0=x_{i}s_{i+1}x_{i}$ and
$v_1=(m-1)(s'_{i+2},s_{i+1})$, where $s'_{i+2}=s'_{i+1}<s_{i+1}$, a
contradiction. Then, $|u_i|$=1 and $v_0=u_iu_{i+1}\cdots u_j$ such
that $|u_i|=\cdots=|u_j|=1$ for some $j$. Then $v_2\cdots
v_{q+1}=u_{j+1}\cdots u_l$ and by Lemma \ref{ls.2},
$|u_{j+1}|=\cdots =|u_l|=1$. Moreover, $|u_{l+1}|=1$.

Case 3. $v_0=(m-1)(x_{i},s_{i+1})$ and $|v_0|=2$, i.e.,
$v_0=x_{i}s_{i+1}$. If $|u_{i+1}|>1$, we have
$v_1=(m-1)(s'_{i+1},x_{i})$. If $i=0$, then $x_{0}=s_{0}$ and $
s'_1=s'_0$. If $|v_1|>1$, then $s_{2}=x_{0}=s_{0}$, a contradiction.
Then $|u_0|=|v_1|=1, p_0=u_0$, a contradiction. Then $i>0$.
Moreover, $s'_{i+1}=s_{i}$, $m_{s_{i}s'_{i}}$ is odd,
$x_{i}=s_{i+2}$ and $u_{i+1}=s_{i+1}s'_{i+1}$. Then
$s'_{i+2}=s_{i+1}$ and $s'_{i+1}<s_{i+2}<s'_{i+2}=s_{i+1}$, also a
contradiction. Thus $|u_{i+1}|=1$ and $u_{i+2}\cdots u_l=v_1\cdots
v_{q+1}$. By Lemma \ref{ls.2}, we have $|u_{i+2}|=\cdots=|u_l|=1$.
\hfill $\blacksquare$

\begin{lemma} \label{s.10}
If $\bar f=u_0\cdots u_i(m-2)(s_i,s'_i)\bar
g(m-2)(s'_{l+1},s_{l+1})u_{l+2}\cdots u_{k+1} y_{k+1}$, and either
$|u_i|= 1$ or $|u_{l+1}|= 1$, then $(f,g)_{\bar f}\equiv 0$.
\end{lemma}
\textbf{Proof} There are two cases to consider.

Case 1. $|u_i|=1$. Suppose $p'_0=s_j$. Then $g=s_i\cdots s_{l+1}-
s_js_i\cdots s_l$.

If $j=l+1$, i.e., $p'_0=s_{l+1}$, then we have $s'_{l+1}\rhd s_i\rhd
s_{l+1},\ s_{l+1}\rhd s_n,\ s'_{l+1}\rhd s_n$ for all $n,\ i+1\leq
n\leq l$, and
\begin{eqnarray*}
(f,g)_{\bar f}&\equiv& u_0\cdots u_{i-1}s_{l+1}s_i\cdots
s_l(m-2)(s'_{l+1},s_{l+1})\cdots u_{k+1}y_{k+1}\\
&&-s'_0u_0\cdots u_{i-1}s_{l+1}s_i\cdots
s_l(m-2)(s'_{l+1},s_{l+1})\cdots u_{k+1}\\
&\equiv& u_0\cdots u_{i-1}s_{l+1}s'_{l+1}s_i\cdots
s_l(m-3)(s_{l+1},s'_{l+1})\cdots u_{k+1}y_{k+1}\\
&&-s'_0u_0\cdots u_{i-1}s_{l+1}s'_{l+1}s_i\cdots
s_l(m-3)(s_{l+1},s'_{l+1})\cdots u_{k+1}\\
&\equiv& u_0\cdots u_{i-1}u_{l+1}s_i\cdots
s_lu_{l+2}\cdots u_{k+1}y_{k+1}\\
&&-s'_0u_0\cdots u_{i-1}u_{l+1}s_i\cdots
s_lu_{l+2}\cdots u_{k+1}\\
&\equiv& 0
\end{eqnarray*}
since $u_0\cdots u_{i-1}u_{l+1}u_i\cdots u_{l}u_{l+2}\cdots
u_{k+1}y_{k+1}-s_0u_0\cdots u_{i-1}u_{l+1}u_i\cdots
u_{l}u_{l+2}\cdots u_{k+1}$ is in (\ref{e3}).

If $j<l+1$, then there exists $i'$ such that $|v_{i'}|>1$ which
implies $|u_{l+1}|=1$. Then by Lemma \ref{s.6}, $(f,g)_{\bar
f}\equiv 0$.

Case 2. $|u_{l+1}|=1$ and $|u_i|\neq 1$. Suppose $p'_0=s_j$. Then
$x_i\rhd s_j\rhd s_{i+1},\cdots,s_{j-1}$, $s'_{i+1}=s'_{j+1}\rhd
s_j$ and
\begin{eqnarray*}
(f,g)_{\bar f}&\equiv& u_0\cdots
u_{i-1}(m-2)(s_i,s'_i)s_{j}x_i\cdots
s_lu_{l+2}\cdots u_{k+1}y_{k+1}\\
&&-s'_0u_0\cdots u_{i-1}(m-2)(s_i,s'_i)s_{j}x_i\cdots
s_lu_{l+2}\cdots u_{k+1}.
\end{eqnarray*}
Since $(m-2)(s_i,s'_i)s_{j}\mapsto\cdots\mapsto s_j(m-2)(s_i,s'_i)$,
we have
\begin{eqnarray*}
(f,g)_{\bar f}&\equiv& u_0\cdots
u_{i-1}s_{j}(m-1)(s_i,s'_i)s_{i+1}\cdots
s_lu_{l+2}\cdots u_{k+1}y_{k+1}\\
&&-s'_0u_0\cdots u_{i-1}s_{j}(m-1)(s_i,s'_i)s_{i+1}\cdots
s_lu_{l+2}\cdots u_{k+1}\\
&\equiv& 0
\end{eqnarray*}
since $u_0\cdots u_{i-1}u_ju_iu_{i+1}\cdots u_{l}u_{l+2}\cdots
u_{k+1}y_{k+1}-s'_0u_0\cdots u_{i-1}u_{j}u_is_{i+1}\cdots
u_{l}u_{l+1}\cdots u_{k+1}$ is in (\ref{e3}). \hfill $\blacksquare$

\ \

\begin{theorem}  \label{s.11}
Suppose that $\bar f=u_0\cdots u_i(m-2)(s_i,s'_i)\bar
g(m-2)(s'_{l+1},s_{l+1})u_{l+2}\cdots u_{k+1} y_{k+1}$,  $|u_i|> 1$
and  $|u_{l+1}|>1$. Then one of the following holds:
\begin{enumerate}
\item[(i)]  $|v_n|=1$ for all $n,\ 0\leq n\leq q+1$
 and
\begin{eqnarray*}
(f,g)_{\bar f}&\equiv& u_0\cdots
u_{i-1}(m-2)(s_i,s'_i)s_{l+1}x_is_{i+1}\cdots
s_l(m-2)(s'_{l+1},s_{l+1})u_{l+2}\cdots u_{k+1}y_{k+1}\\
&&-s'_0u_0\cdots u_{k+1}.
\end{eqnarray*}

\item[(ii)] $|v_0|=2,\ |v_n|=1$ for all $n,\ 1\leq n\leq q+1$ and

if $m_{s_is'_i}=3$, then $(f,g)_{\bar f}\equiv 0$;

if $m_{s_is'_i}>3$, then $s_{l+1}=x_i$, $s'_{l+1}=s'_{i+1}=y_i$ and
\begin{eqnarray*}
(f,g)_{\bar f}&\equiv& u_0\cdots
u_{i-1}(m-3)(s_i,s'_i)s_{i+1}s'_{i+1}x_is_{i+1}\cdots
s_l(m-2)(s'_{l+1},s_{l+1})u_{l+2}\cdots u_{k+1} y_{k+1}\\
&&-s'_0u_0\cdots u_{k+1}.
\end{eqnarray*}
\end{enumerate}
\end{theorem}
\textbf{Proof} By Lemma \ref{s.9}, we have $|u_{i+1}|=\cdots
=|u_{l}|=1$ and so $\bar g=x_is_{i+1}\cdots s_{l}s_{l+1}$. There are
two cases to consider.

Case 1. $v_0=x_i$. Then $|v_n|=1$ for all $n,\ 1\leq n\leq q+1$.
Otherwise, $z_{q+1}=s_{l+1}\in \{s_{i+1},\cdots, s_l\}$ which shows
$|u_{l+1}|=1$, a contradiction. Then
\begin{eqnarray*}
(f,g)_{\bar f}&\equiv& u_0\cdots
u_{i-1}(m-2)(s_i,s'_i)s_{l+1}x_is_{i+1}\cdots
s_l(m-2)(s'_{l+1},s_{l+1})u_{l+2}\cdots u_{k+1} y_{k+1}\\
&&-s'_0u_0\cdots u_{k+1}.
\end{eqnarray*}

Case 2. $v_0=x_is_{i+1}$. Then $|v_n|=1$ for all $n,\ 1\leq n\leq
q+1$. Otherwise, we have $s_j=p_0=x_i$ for some $j\ \ (i+1<j<l+1)$.
Then $u_j=(m-1)(x_i,s'_{i+1})$ and $|u_j|=|u_i|>1$, a contradiction.
Therefore $z_{q+1}=x_i=s_{l+1}$, $u_{l+1}=(m-1)(x_i,s'_{i+1})$ and
\begin{eqnarray*}
(f,g)_{\bar f}&\equiv& u_0\cdots
u_{i-1}(m-2)(s_i,s'_i)s_{i+1}x_is_{i+1}\cdots
s_l(m-2)(s'_{l+1},s_{l+1})u_{l+2}\cdots u_{k+1}y_{k+1}\\
&&-s'_0u_0\cdots u_{k+1}\\
 &\equiv&
u_0\cdots u_{i-1}(m-3)(s_i,s'_i)s_{i+1}s'_{i+1}x_is_{i+1}\cdots
s_l(m-2)(s'_{l+1},s_{l+1})u_{l+2}\cdots u_{k+1} y_{k+1}\\
&&-s'_0u_0\cdots u_{k+1}.
\end{eqnarray*}
If $m_{s_is'_i}=3$, we have $s'_{i+1}=s_i,\
x_i=s'_{i}=s_{l+1}<s'_{l+1}=s'_{i+1}=s_i$. Therefore $i=0$ and
\begin{eqnarray*}
(f,g)_{\bar f}&\equiv&  s_{1}s_{0}s'_0s_{1}\cdots
s_ls_0u_{l+2}\cdots u_{k+1} y_{k+1}-s'_0s_{1}s_{0}s'_0s_{1}\cdots
s_ls_0u_{l+2}\cdots u_{k+1}\\
&\equiv&s_{1}s'_0s_{0}s'_0s_{1}\cdots s_lu_{l+2}\cdots u_{k+1}
y_{k+1}-s'_0s_{1}s'_0s_{0}s'_0s_{1}\cdots s_lu_{l+2}\cdots
u_{k+1}\\
&\equiv&s_{1}s'_0s_{0}s_1s'_0s_{1}\cdots s_lu_{l+2}\cdots
u_{k+1}-s_{1}s'_0s_1s_{0}s'_0s_{1}\cdots s_lu_{l+2}\cdots
u_{k+1}\\
&\equiv&s_{1}s'_0s_1s_{0}s'_0s_{1}\cdots s_lu_{l+2}\cdots
u_{k+1}-s_{1}s'_0s_1s_{0}s'_0s_{1}\cdots s_lu_{l+2}\cdots
u_{k+1}\\
&\equiv& 0.
\end{eqnarray*}

If $m_{s_is'_i}>3$, we have
\begin{eqnarray*}
(f,g)_{\bar f}&\equiv&  u_0\cdots
u_{i-1}(m-3)(s_i,s'_i)s_{i+1}s'_{i+1}x_is_{i+1}\cdots
s_l(m-2)(s'_{l+1},s_{l+1})u_{l+2}\cdots u_{k+1} y_{k+1}\\
&&-s'_0u_0\cdots u_{k+1}.
\end{eqnarray*}
The proof is completed. \hfill $\blacksquare$

\ \

 \begin{lemma}  \label{s.12}
Suppose $\bar f=u_0\cdots u_{i-1}(m-3)(s_i,s'_i)\bar
g(m-2)(s'_{l+1},s_{l+1})u_{l+2}\cdots u_{k+1} y_{k+1}$. Then
$|u_{i+1}|=\cdots =|u_l|=1$.
\end{lemma}
\textbf{Proof} Clearly, $v_0=s'_{i+1}$. There are two cases to
consider.

Case 1. $v_1=x_{i}$. Since $u_{i+1}\cdots u_l=v_2\cdots v_{q+1}$, by
Lemma \ref{ls.2}, we have $|u_{i+1}|=\cdots=|u_l|=1$.

Case 2. $v_1=(m-1)(x_{i},s_{i+1})$. We have $m_{s_{i+1}s'_{i+1}}=2$,
i.e., $|u_{i+1}|=1$.

If $|v_1|>2$, then $s_{i+2}=x_{i}$. We have $|u_{i+2}|>1$,
$|v_1|=3$, $v_2=(m-1)(s'_{i+2},s_{i+1})$ and
$s'_{i+1}=s'_{i+2}<s_{i+1}$, a contradiction. Then $|v_1|=2$ and
$v_2\cdots v_{q+1}=u_{i+2}\cdots u_l$. By Lemma \ref{ls.2}, we have
$|u_{i+2}|=\cdots =|u_l|=1$. \hfill $\blacksquare$

\ \

\begin{theorem}  \label{s.14}
Suppose that $\bar f=u_0\cdots u_i(m-3)(s_i,s'_i)\bar
g(m-2)(s'_{l+1},s_{l+1})u_{l+2}\cdots u_{k+1} y_{k+1}$ and
$|u_{l+1}|>1$. Then one of the following holds.
\begin{enumerate}
\item[(i)]  $|v_0|=|v_1|=1$ and

if $i=0$, then $(f,g)_{\bar f}\equiv 0$;

if $i>0$, then
\begin{eqnarray*}
(f,g)_{\bar f}&\equiv& u_0\cdots
u_{i-1}s_is_j(m-2)(s'_{i},s_{i})s_{i+1}\cdots s_{l}u_{l+2}\cdots
u_{k+1}y_{k+1}-s'_0u_0\cdots  u_{k+1}.
\end{eqnarray*}

\item[(ii)]  $|v_0|=1,\ |v_1|=2,\ |v_n|=1$ for all $n \ (1<n\leq q+1)$ and
\begin{eqnarray*}
(f,g)_{\bar f}&\equiv& u_0\cdots
u_{i-1}(m-3)(s_{i},s'_{i})s_{i+1}s'_{i+1}x_{i}s_{i+1}\cdots
s_{l}(m-2)(s'_{i+1}x_{i})u_{l+2}\cdots u_{k+1}y_{k+1}\\
&&-s'_0u_0\cdots u_{k+1}.
\end{eqnarray*}

\end{enumerate}
\end{theorem}
\textbf{Proof} By Lemma \ref{s.12}, $|u_{i+1}|=\cdots=|u_{l}|=1$ and
$v_0=s'_{i+1}$. There are two cases to consider.

Case 1. $v_1=x_i$. There exists $s_j=p'_0$, where $s_j\in
\{s_{i+1},\cdots, s_{l+1}\}$. If $j=l+1$, then $|v_{n}|=1$ for all
$n$ and $s'_{l+1}=s'_{i+1}\rhd s_{l+1}$. Thus, $|u_{l+1}|=1$. If
$j\neq l+1$, there exists $|v_{n}|>1\ (n>1)$. Then
$s_{l+1}\in\{s_{i+1},\cdots,s_l\}$ and $|u_{l+1}|=1$.

If $i>0$, then $s'_{i+1}=s'_{i},\ x_{i}=s_{i}$. Hence
$m_{s_{i}s'_{i}}$ is even and $s'_{i+1}=\cdots =s'_{l+1}=s'_{i}$.
Now by ELW's, we have
\begin{eqnarray*}
(f,g)_{\bar f}&\equiv& u_0\cdots
u_{i-1}(m-3)(s_{i},s'_{i})s_js'_{i}s_{i}s_{i+1}\cdots
s_{l}u_{l+2}\cdots u_{k+1}y_{k+1}-s'_0u_0\cdots u_{k+1}\\
 &\equiv& u_0\cdots
u_{i-1}s_is_j(m-2)(s'_{i},s_{i})s_{i+1}\cdots s_{l}u_{l+2}\cdots
u_{k+1}y_{k+1}-s'_0u_0\cdots u_{k+1}.
\end{eqnarray*}

If $i=0$,  $s_1'=s_0$ since $s'_{1}>x_0$. Then $m_{s_0s'_0}$ is odd
and $s_0=s'_{1}=s'_2=\cdots s'_{l+1}=s'_{l+2}$. Since
$h=u_0s_1\cdots s_{l}u_{l+2}\cdots u_{k+1}y_{k+1}-s'_0u_0s_1\cdots
s_{l}u_{l+2}\cdots u_{k+1}$ is in (\ref{e3}), we have
\begin{eqnarray*}
(f,g)_{\bar f}&\equiv& (m-3)(s_0,s'_0)s_{j}s_0s'_0s_1\cdots
s_{l}u_{l+2}\cdots
u_{k+1}y_{k+1}\\
&& -(m-2)(s'_0,s_0)s_{j}s_0s'_0s_1\cdots s_{l}u_{l+2}\cdots
u_{k+1}\\
&\equiv& s_{j}u_0s_1\cdots s_{l}u_{l+2}\cdots
u_{k+1}y_{k+1}-s'_0s_ju_0s_1\cdots s_{l}u_{l+2}\cdots
u_{k+1}\\
&\equiv& s_{j}u_0s_1\cdots s_{l}u_{l+2}\cdots
u_{k+1}y_{k+1}-s_js'_0u_0s_1\cdots s_{l}u_{l+2}\cdots u_{k+1}\\
&\equiv& 0.
\end{eqnarray*}

Case 2. $v_1=x_is_{i+1}$. Clearly,
$x_{i}\not\in\{s_{i+1},\cdots,s_{l}\}$. Otherwise, $x_{i}=s_{j}$ for
some $j\ (i+1\leq j\leq l)$ and so $|u_{j}|=|u_i|>1$, a
contradiction. Then $x_{i}=s_{l+1}$, $|v_2|=\cdots=|v_{q+1}|=1$ and
$u_{l+1}=(m-1)(x_i,s'_{i+1})$. Moreover, we have $s'_{i+1}\rhd
s_{i+1}>x_i$ and $m_{s_{i}s'_{i}}>2,\ m_{x_{i}s_{i+1}}=3$. By ELW's,
we have
\begin{eqnarray*}
(f,g)_{\bar f}&\equiv&u_0\cdots
u_{i-1}(m-3)(s_{i},s'_{i})s_{i+1}s'_{i+1}x_{i}s_{i+1}\cdots
s_{l}(m-2)(s'_{i+1}x_{i})u_{l+2}\cdots u_{k+1}y_{k+1}\\
&&-s'_0u_0\cdots u_{k+1}.
\end{eqnarray*}
If $m_{s_is'_i}=3$, then $|u_i|=2$, $x_i=s'_i$ and $s'_{i+1}=s_{i}$.
Since $s'_{i+1}>s_{i+1}>x_i$, we have $i=0$ and $a=1$, which
contradicts $a\neq 1$. Then $|u_i|>2$ and $m_{s_is'_i}$ is even.
\hfill $\blacksquare$

 \ \

Now, let
\begin{eqnarray*}\label{e8}
g=m(s,s')-m(s',s),\ \ s>s'
\end{eqnarray*}
be a relation in (\ref{e2}) and $f$ be as (\ref{e6}) again. In the
following Lemma (Theorems) \ref{s.15}--\ref{s.17}, we will deal with
another inclusion compositions $(f,g)_w,\ w=\bar f=a\bar gb,\ \bar
g=m(s,s')$. There are another two nontrivial cases which will be
mentioned in Theorems \ref{s.16} and \ref{s.17}.

\begin{lemma}  \label{s.15}
If $\bar f=u_0\cdots u_{i-1}\bar gu_{l+2}\cdots u_{k+1} y_{k+1}$,
then $|u_{i}|\cdots=|u_{l+1}|=1$ and $(f,g)_{\bar f}\equiv 0$.
\end{lemma}
\textbf{Proof} If there exists $j$ such that $|u_j|>1$, there will
be three different letters in $\bar g$, a contradiction. Therefore,
$|u_i|=\cdots =|u_{l+1}|=1$.

Similarly to the proof of Lemma \ref{s.6}, the result holds. \hfill
$\blacksquare$

\ \

\begin{theorem}  \label{s.16}
Suppose $\bar f=u_0\cdots u_{i-1}(m-2)(s_i,s'_i)\bar g
(m-2)(s'_{i+1},s_{i+1})u_{i+2}\cdots u_{k+1} y_{k+1}$, where $0\leq
i\leq k$, $u_0\cdots u_{i-1}=1$ if $i=0$, $u_{i+2}\cdots u_{k+1}=1$
if $i\leq k$. Then the following statements hold.
\begin{enumerate}
\item[(i)]\ $g=x_is_{i+1}-s_{i+1}x_i,\
x_i>s_{i+1}$.
\item[(ii)]\ $(f,g)_{\bar f}\equiv 0$ if $|u_i|=1$ or
$|u_{i+1}|=1$.
\item[(iii)]\ If $|u_i|>1$ and $|u_{i+1}|>1$, then
\begin{eqnarray*}
(f,g)_{\bar f}&\equiv&u_0\cdots u_{i-1}(m-2)(s_i,s'_i)s_{i+1}x_i
(m-2)(s'_{i+1},s_{i+1})u_{i+2}\cdots u_{k+1} y_{k+1}\\
&&-s'_0u_0\cdots  u_{k+1}.
\end{eqnarray*}
\end{enumerate}
\end{theorem}
\textbf{Proof} (i) is clear.

(ii) Suppose $|u_i|=1$. Since $u_0\cdots
u_{i-1}u_{i+1}s_iy_{i+1}-s'_0u_0\cdots u_{i-1}u_{i+1}s_i$ is in
(\ref{e3}), we have
\begin{eqnarray*}
(f,g)_{\bar f}&\equiv& u_0\cdots u_{i-1}s_{i+1}s_i
(m-2)(s'_{i+1},s_{i+1})u_{i+2}\cdots u_{k+1} y_{k+1}\\
&&-s'_0u_0\cdots u_{i-1}s_{i+1}s_i
(m-2)(s'_{i+1},s_{i+1})u_{i+2}\cdots u_{k+1}\\
&\equiv& u_0\cdots u_{i-1}u_{i+1}s_i u_{i+2}\cdots u_{k+1}
y_{k+1}-s'_0u_0\cdots u_{i-1}u_{i+1}s_iu_{i+2}\cdots u_{k+1}\\
&\equiv& u_0\cdots u_{i-1}u_{i+1}s_iy_{i+1} u_{i+2}\cdots u_{k+1}
-s'_0u_0\cdots u_{i-1}u_{i+1}s_iu_{i+2}\cdots u_{k+1}\\
&\equiv& 0.
\end{eqnarray*}

Suppose $|u_{i+1}|=1$.  Since $u_0\cdots
u_{i-1}u_{i+1}s_iy_{i+1}-s'_0u_0\cdots u_{i-1}u_{i+1}s_i$ is in
(\ref{e3}), we have
\begin{eqnarray*}
(f,g)_{\bar f}&\equiv&u_0\cdots u_{i-1}(m-2)(s_i,s'_i)s_{i+1}x_i
u_{i+2}\cdots u_{k+1} y_{k+1}\\
&&-s'_0u_0\cdots u_{i-1}(m-2)(s_i,s'_i)s_{i+1}x_i
u_{i+2}\cdots u_{k+1}\\
&\equiv& u_0\cdots u_{i-1}u_{i+1}s_i u_{i+2}\cdots u_{k+1}
y_{k+1}-s'_0u_0\cdots u_{i-1}u_{i+1}s_iu_{i+2}\cdots u_{k+1}\\
&\equiv& u_0\cdots u_{i-1}u_{i+1}s_iy_{i+1} u_{i+2}\cdots u_{k+1}
-s'_0u_0\cdots u_{i-1}u_{i+1}s_iu_{i+2}\cdots u_{k+1}\\
&\equiv& 0.
\end{eqnarray*}

(iii) Suppose $|u_i|>1$ and $|u_{i+1}|>1$. Then by ELW's, we have
\begin{eqnarray*}
(f,g)_{\bar f}\equiv u_0\cdots u_{i-1}(m-2)(s_i,s'_i)s_{i+1}x_i
(m-2)(s'_{i+1},s_{i+1})u_{i+2}\cdots u_{k+1} y_{k+1}-s'_0u_0\cdots
u_{k+1}.
\end{eqnarray*}
The proof is completed. \hfill $\blacksquare$

\ \

\begin{theorem}  \label{s.17}
Suppose $\bar f=u_0\cdots u_{i-1}(m-2)(s_i,s'_i)\bar g
(m-2)(s'_{i+2},s_{i+2})u_{i+3}\cdots u_{k+1} y_{k+1}$, where $0\leq
i\leq k-1$, $u_0\cdots u_{i-1}=1$ if $i=0$, $u_{i+3}\cdots
u_{k+1}=1$ if $i=k-1$. Then the following statements hold.
\begin{enumerate}
\item[(i)]\ $g=x_is_{i+1}x_i-s_{i+1}x_is_{i+1},\
x_i>s_{i+1}$.
\item[(ii)]\ $(f,g)_{\bar f}\equiv 0$  if $|u_i|=1$
or $|u_{i}|=2$.
\item[(iii)]\ If $|u_i|>2$, then
\begin{eqnarray*}
(f,g)_{\bar f}&\equiv &u_0\cdots
u_{i-1}(m-3)(s_i,s'_i)s_{i+1}s'_{i+1}x_is_{i+1}
(m-2)(s'_{i+2},x_i)u_{i+3}\cdots u_{k+1} y_{k+1}\\
&&-s'_0u_0\cdots  u_{k+1}.
\end{eqnarray*}
\end{enumerate}
\end{theorem}
\textbf{Proof} (i) is clear.

(ii) If $|u_i|=1$, then $|u_{i+1}|=|u_{i+2}|=1$. Similar to Lemma
\ref{s.6}, we have $(f,g)_{\bar f}\equiv 0$.

If $|u_i|=2$, then $u_i=s_is'_i$,
$s'_i=x_i=s_{i+2}<s'_{i+2}=s'_{i+1}=s_i$ which implies $i=0$. Since
$s_0s'_0s_1s_0-s'_0s_0s'_0s_1$ and $s'_0s_1u_3\cdots
u_{k+1}y_{k+1}-s_1s'_0s_1u_3\cdots u_{k+1}$ are in (3),
\begin{eqnarray*}
(f,g)_{\bar f}&\equiv& s_0s_1s'_0s_1s_0u_3\cdots
u_{k+1}y_{k+1}-s'_0s_0s_1s'_0s_1s_0u_3\cdots u_{k+1}\\
&\equiv&s_1s_0s'_0s_1s_0u_3\cdots
u_{k+1}y_{k+1}-s'_0s_1s_0s'_0s_1s_0u_3\cdots u_{k+1}\\
&\equiv&s_1s'_0s_0s'_0s_1u_3\cdots
u_{k+1}y_{k+1}-s'_0s_1s'_0s_0s'_0s_1u_3\cdots u_{k+1}\\
&\equiv&s_1s'_0s_0s_1s'_0s_1u_3\cdots
u_{k+1}-s_1s'_0s_1s_0s'_0s_1u_3\cdots u_{k+1}\\
&\equiv&s_1s'_0s_1s_0s'_0s_1u_3\cdots
u_{k+1}-s_1s'_0s_1s_0s'_0s_1u_3\cdots u_{k+1}\\
 &\equiv& 0.
\end{eqnarray*}

(iii) Suppose $|u_i|>2$. By ELW's, we have
\begin{eqnarray*}
(f,g)_{\bar f}&\equiv& u_0\cdots
u_{i-1}(m-2)(s_i,s'_i)s_{i+1}x_is_{i+1}
(m-2)(s'_{i+2},x_i)u_{i+3}\cdots u_{k+1} y_{k+1}\\
&&-s'_0u_0\cdots  u_{k+1}\\
 &\equiv& u_0\cdots
u_{i-1}(m-3)(s_i,s'_i)s_{i+1}s'_{i+1}x_is_{i+1}
(m-2)(s'_{i+2},x_i)u_{i+3}\cdots u_{k+1} y_{k+1}\\
&&-s'_0 u_0\cdots  u_{k+1}.
\end{eqnarray*}
The proof is completed. \hfill $\blacksquare$

\ \

Now we finish all the cases of inclusion compositions. Most of them
are trivial except six cases which are mentioned in Theorems
\ref{s.11}, \ref{s.14}, \ref{s.16} and \ref{s.17}. But in fact, we
can classify these six cases into four cases.

\ \

Now we consider that in what  instances the nontrivial cases may
happen.

The first nontrivial case, which is the first case of Theorem
\ref{s.11} and the nontrivial case of Theorem \ref{s.16},  happens
if the following $f$ exists:

{\bf C1}: $f=(m-1)(s_0,s'_0)\cdots (m-1)(s_i,s'_i)\cdots
(m-1)(s_{l+1},s'_{l+1})\cdots m(s_{k+1},s'_{k+1})-m(s'_0,s_0)\cdots
(m-1)(s_{k+1},s'_{k+1})$, where $0\leq i\leq l\leq k$, such that
\begin{eqnarray*}
&&(a)\
|(m-1)(s_i,s'_i)|\geq 2, \ |(m-1)(s_{l+1},s'_{l+1})|\geq 2,\ x_i\rhd s_{l+1};\\
&&(b)\ (m-1)(s_j,s'_j)=s_j, \ \ s_{l+1}\rhd s_j\ \mbox{ for any }j,\
i+1\leq j\leq l.
\end{eqnarray*}
\noindent {\bf Remarks}: In the case {\bf C1}, we have

1) $\bar f$ contains $\bar g$ as a subword where $g=x_is_{i+1}\cdots
s_{l+1}-s_{l+1}s_{i+1}\cdots s_l$, $x_i\rhd s_{l+1}$ and
$s_{l+1}\rhd s_j\ \mbox{ for any }j,\ i+1\leq j\leq l$.

2) If there is no  $f\in (3')$ with {\bf C1} where $0\leq i= l\leq
k$, then for any $f\in (3'), \ f$ is not with property {\bf C1}.

\ \

The second nontrivial case, which is the second case of Theorem
\ref{s.11} and the nontrivial case of Theorem \ref{s.17}, happens if
the following $f$ exists:

{\bf C2}: $f=(m-1)(s_0,s'_0)\cdots (m-1)(s_i,s'_i)\cdots
(m-1)(s_{l+1},s'_{l+1})\cdots m(s_{k+1},s'_{k+1})-m(s'_0,s_0)\cdots
(m-1)(s_{k+1},s'_{k+1})$, where $0\leq i<l\leq k$, such that
\begin{eqnarray*}
&&(a)\ |(m-1)(s_i,s'_i)|>2,\  (m-1)(s_{i+1},s'_{i+1})=s_{i+1},\ x_i=s_{l+1}>s_{i+1},\ m_{x_is_{i+1}}=3;\\
&&(b)\ (m-1)(s_j,s'_j)=s_j,\  s_{l+1}\rhd s_j\mbox{ for any }j,\
i+2\leq j\leq l.
\end{eqnarray*}
\noindent {\bf Remarks}: In the case {\bf C2}, we have

1) $\bar f$ contains $\bar g$ as a subword where $g=x_is_{i+1}\cdots
s_{l+1}-s_{l+1}s_{i+1}\cdots s_l$, $m_{x_is_{i+1}}=3$ and
$s_{l+1}\rhd s_j\mbox{ for any }j,\ i+2\leq j\leq l$.

2) If there is no  $f\in (3')$ with {\bf C2}  where  $0\leq i=
l-1\leq k-1$, then for any $f\in (3'), \ f$ is not with property
{\bf C2}.

\ \

The third nontrivial case, which is the first case of Theorem
\ref{s.14}, happens if the following $f$ exists:

 {\bf C3}:
$f=(m-1)(s_0,s'_0)\cdots (m-1)(s_i,s'_i)\cdots
(m-1)(s_{l+1},s'_{l+1})\cdots m(s_{k+1},s'_{k+1})-m(s'_0,s_0)\cdots
(m-1)(s_{k+1},s'_{k+1})$, where $0\leq i\leq l\leq k$, such that
\begin{eqnarray*}
&&(a)\ (m-1)(s_{i},s'_{i})\geq 2,\ m_{s_is'_i}\ \mbox{is\ even and
there\ exists\ }m \ (\ i+1\leq m\leq l+1)\ \mbox{such \
that\ }\\
&&\ \ \ \ s'_{i+1}\rhd s_m\rhd x_i;\\
&&(b)\ (m-1)(s_j,s'_j)=s_j\mbox{ for any } j,\ i+1\leq j\leq l+1,\\
&&\ \ \ \ s_m\rhd s_n \mbox{ for any }n,\  i+1 \leq n\leq
m-2\ \mbox{and}\\
&&\ \ \ \ s_{m-1}\cdots s_{i_1}=(m-1)(s_{m-1},s_m),\ s_{m-1}<s_m,\\
&&\ \ \ \ s_{i_1+1}\cdots
s_{i_2}=(m-1)(s_{i_1+1},s_{i_1+2}),\ s_{i_1+1}<s_{i_1+2},\cdots,\\
&&\ \ \ \  s_{i_n+1}\cdots s_{l+1}=m(s_{i_n+1},s_{i_n+2}),\
s_{i_n+1}<s_{i_n+2}.
\end{eqnarray*}
\noindent {\bf Remarks}: In the case {\bf C3}, we have

1) $\bar f$ contains $\bar g$ as a subword where
$g=s'_{i+1}x_is_{i+1}\cdots s_{l+1}-s_ms'_{i+1}x_is_{i+1}\cdots
s_{l}\in(3')$ such that $s'_{i+1}\rhd s_m\rhd x_i$ for some $m\ (\
i+1\leq m\leq l+1)$ and $s_m\rhd s_n$ for any $n,\  i+1 \leq n\leq
m-2.$

2) If there is no  $f\in (3')$ with {\bf C3}  where  $0\leq i= l\leq
k$, then for any $f\in (3'), \ f$ is not with property {\bf C3}.

\ \

The fourth nontrivial case, which is the second case of Theorem
\ref{s.14}, happens if the following $f$ exists:

{\bf C4}: $f=(m-1)(s_0,s'_0)\cdots (m-1)(s_i,s'_i)\cdots
(m-1)(s_{l+1},s'_{l+1})\cdots m(s_{k+1},s'_{k+1})-m(s'_0,s_0)\cdots
(m-1)(s_{k+1},s'_{k+1})$, where $0\leq i<l\leq k$, such that
\begin{eqnarray*}
&&(a)\ |(m-1)(s_i,s'_i)|>2,\ (m-1)(s_{i+1},s'_{i+1})=s_{i+1},\ s_{i+1}>x_i=s_{l+1},\ m_{x_is_{i+1}}=3;\\
&&(b)\ (m-1)(s_j,s'_j)=s_j,\ s_{l+1}\rhd s_j \mbox{ for any }j,\
i+2\leq j\leq l.
\end{eqnarray*}
\noindent {\bf Remarks}: In the case {\bf C4}, we have

1) $\bar f$ contains $\bar g$ as a subword where
$g=s'_{i+1}x_i\cdots s_{l+1}-s_{i+1}s'_{i+1}x_i\cdots s_l\in (3')$,
$s'_{i+1}\rhd s_{i+1},\ m_{x_is_{i+1}}=3$ and $s_{l+1}\rhd s_j$ for
any $j,\ i+2\leq j\leq l$.

2) If there is no  $f\in (3')$ with {\bf C4}  where  $0\leq i=
l-1\leq k-1$, then for any $f\in (3'), \ f$ is not with property
{\bf C4}.

\ \

\noindent\textbf{Remark:} In the Example \ref{ex1}, there exist
relations in $(3')$ with properties {\bf C1} and {\bf C2}.

\begin{theorem}\label{t3.19}
$S=\{(\ref{e1}),(\ref{e2}),(3')\}$ is a Gr\"{o}bner-Shirshov basis
of $W$ if there is no $f\in (3')$ with properties ${\bf C1}\vee {\bf
C2}\vee {\bf C3}\vee {\bf C4}$.
\end{theorem}

\textbf{Proof}
 We will prove that all possible compositions are trivial
modulo $S$. Denote by $(i\wedge j)_w$ the composition of the type
$(i)$ and type $(j)$ with respect to the ambiguity $w$.

By Lemmas \ref{s.6}, \ref{s.10} and \ref{s.15}, and Theorems
\ref{s.11}, \ref{s.14}, \ref{s.16} and \ref{s.17}, we know that all
inclusion compositions are trivial. Thus, we need only to check the
intersection compositions.

\begin{enumerate}
\item[($1\wedge2$)]\ $w=sm(s,s'),\ s>s'$.
\begin{eqnarray*}
(1\wedge2)_w&=&-(m-1)(s',s)+sm(s',s)\\
&\equiv&-(m-1)(s',s)+(m+1)(s,s')\\
&\equiv&-(m-1)(s',s)+(m-1)(s',s)\\
&\equiv&0.
\end{eqnarray*}

\item[($1\wedge3'$)]\
$w=s_0(m-1)(s_0,s'_0)(m-1)(s_1,s'_1)\cdots
(m-1)(s_k,s'_k)m(s_{k+1},s'_{k+1}).$
\begin{eqnarray*}
(1\wedge3)_w&=&-(m-2)(s'_0,s_0)(m-1)(s_1,s'_1)\cdots(m-1)(s_{k},s'_{k})m(s_{k+1},s'_{k+1})\\
&&+s_0m(s'_0,s_0)(m-1)(s_1,s'_1)\cdots
(m-1)(s_{k},s'_{k})(m-1)(s_{k+1},s'_{k+1})\\
&\equiv&-(m-1)(s'_0,s_0)(m-1)(s_1,s'_1)\cdots
(m-1)(s_{k},s'_{k})(m-1)(s_{k+1},s'_{k+1})\\
&&+(m-1)(s'_0,s_0)(m-1)(s_1,s'_1)\cdots
(m-1)(s_{k},s'_{k})(m-1)(s_{k+1},s'_{k+1})\\
&\equiv&0.
\end{eqnarray*}

\item[($2\wedge1$)]\ $w=m(s,s')x,\ s>s'$, where $x$ is the last letter of $m(s,s')$.
\begin{eqnarray*}
(2\wedge1)_w&=&-m(s',s)x+(m-1)(s,s')\\
&\equiv&-(m+1)(s',s)+(m-1)(s,s')\\
&\equiv&-(m-1)(s,s')+(m-1)(s,s')\\
&\equiv&0.
\end{eqnarray*}

\item[($2\wedge2$)]\ There are two cases to consider.

Case 1. $w=m(s,s')(m-1)(s'',x),\ s>s',\ x>s''$, where $x$ is the
last letter of $m(s,s')$.
$$
(2\wedge2)_w=-m(s',s)(m-1)(s'',x)+(m-1)(s,s')m(s'',x)\equiv0.
$$

Case 2. $w=(2i)(s,s')m(s,s'),\  s>s',\ 1\leq i< m_{ss'}/2$. We just
prove the case that $m_{ss'}$ is even. For the case that $m_{ss'}$
is odd, the proof is similar. Assume that $m_{ss'}$ is even. Then
$$
(2\wedge2)_w\equiv
-m(s',s)(2i)(s,s')+(2i)(s,s')m(s',s)\equiv-(m-2i)(s',s)+(m-2i)(s',s)\equiv
0.
$$

\item[($2\wedge3'$)]\ There are two cases to consider.

Case 1. $ w=(m-1)(s,s')(m-1)(s_0,s'_0)(m-1)(s_1,s'_1)\cdots
(m-1)(s_k,s'_k)m(s_{k+1},s'_{k+1})$, $s>s',\ s_0>s'_0$, where $s_0$
is the last letter of $m(s,s')$. Since
$h=(m-1)(s,s')m(s'_0,s_0)-m(s',s)(m-1)(s'_0,s_0)\in (3')$, we have
\begin{eqnarray*}
&&(2\wedge3')_w\\
&=&-m(s',s)(m-2)(s'_0,s_0)(m-1)(s_1,s'_1)\cdots
(m-1)(s_k,s'_k)m(s_{k+1},s'_{k+1})\\
&&+(m-1)(s,s')m(s'_0,s_0)(m-1)(s_1,s'_1)\cdots
(m-1)(s_{k},s'_{k})(m-1)(s_{k+1},s'_{k+1})\\
&\equiv&-m(s',s)(m-1)(s'_0,s_0)(m-1)(s_1,s'_1)\cdots
(m-1)(s_{k},s'_{k})(m-1)(s_{k+1},s'_{k+1})\\
&&+m(s',s)(m-1)(s'_0,s_0)(m-1)(s_1,s'_1)\cdots
(m-1)(s_{k},s'_{k})(m-1)(s_{k+1},s'_{k+1})\\
&\equiv&0.
\end{eqnarray*}

Case 2. $ w=(2i)(s_0,s'_0)(m-1)(s_0,s'_0)(m-1)(s_1,s'_1)\cdots
(m-1)(s_k,s'_k)m(s_{k+1},s'_{k+1})$, \ $1\leq i<m_{s_0s'_0}/2$. We
prove only the case that $m_{s_0s'_0}$ is even.  For the case that
$m_{s_0s'_0}$ is odd, the proof is similar. Assume that
$m_{s_0s'_0}$ is even. Then
\begin{eqnarray*}
&&(2\wedge3')_w\\
&=&-m(s'_0,s_0)(2i-1)(s_0,s'_0)(m-1)(s_1,s'_1)\cdots(m-1)(s_k,s'_k)m(s_{k+1},s'_{k+1})\\
&&+(2i)(s_0,s'_0)m(s'_0,s_0)(m-1)(s_1,s'_1)\cdots
(m-1)(s_{k},s'_{k})(m-1)(s_{k+1},s'_{k+1})\\
&=&-(m-2i+1)(s'_0,s_0)(m-1)(s_1,s'_1)\cdots(m-1)(s_k,s'_k)m(s_{k+1},s'_{k+1})\\
&&+(m-2i)(s'_0,s_0)(m-1)(s_1,s'_1)\cdots
(m-1)(s_{k},s'_{k})(m-1)(s_{k+1},s'_{k+1})\\
&\equiv&-(m-2i+1)(s'_0,s_0)s'_1(m-1)(s_1,s'_1)\cdots
(m-1)(s_{k},s'_{k})(m-1)(s_{k+1},s'_{k+1})\\
&&+(m-2i)(s'_0,s_0)(m-1)(s_1,s'_1)\cdots
(m-1)(s_{k},s'_{k})(m-1)(s_{k+1},s'_{k+1})\\
&\equiv&-(m-2i)(s'_0,s_0)(m-1)(s_1,s'_1)\cdots
(m-1)(s_{k},s'_{k})(m-1)(s_{k+1},s'_{k+1})\\
&&+(m-2i)(s'_0,s_0)(m-1)(s_1,s'_1)\cdots
(m-1)(s_{k},s'_{k})(m-1)(s_{k+1},s'_{k+1})\\
&\equiv&0.
\end{eqnarray*}

\item[($3'\wedge1$)]
$w=(m-1)(s_0,s'_0)(m-1)(s_1,s'_1)\cdots
(m-1)(s_k,s'_k)m(s_{k+1},s'_{k+1})y_{k+1}$, where $y_{k+1}$ is the
last letter of $m(s_{k+1},s'_{k+1})$.
\begin{eqnarray*}
&&(3'\wedge1)_w\\
&=&-m(s'_0,s_0)(m-1)(s_1,s'_1)\cdots
(m-1)(s_{k},s'_{k})(m-1)(s_{k+1},s'_{k+1})y_{k+1}\\
&&+(m-1)(s_0,s'_0)(m-1)(s_1,s'_1)\cdots
(m-1)(s_{k},s'_{k})(m-1)(s_{k+1},s'_{k+1})\\
&\equiv&-m(s'_0,s_0)s'_1(m-1)(s_1,s'_1)\cdots
(m-1)(s_{k},s'_{k})(m-1)(s_{k+1},s'_{k+1})\\
&&+(m-1)(s_0,s'_0)(m-1)(s_1,s'_1)\cdots
(m-1)(s_{k},s'_{k})(m-1)(s_{k+1},s'_{k+1})\\
&\equiv&-(m-1)(s_0,s'_0)(m-1)(s_1,s'_1)\cdots
(m-1)(s_{k},s'_{k})(m-1)(s_{k+1},s'_{k+1})\\
&&+(m-1)(s_0,s'_0)(m-1)(s_1,s'_1)\cdots
(m-1)(s_{k},s'_{k})(m-1)(s_{k+1},s'_{k+1})\\
&\equiv&0.
\end{eqnarray*}

\item[($3'\wedge2$)]\ There are two cases to consider.

Case 1. $w=(m-1)(s_0,s'_0)(m-1)(s_1,s'_1)\cdots
m(s_{k+1},s'_{k+1})(m-1)(t,y_{k+1}), \ y_{k+1}>t$, where $y_{k+1}$
is the last letter of $m(s_{k+1},s'_{k+1})$. Then
\begin{eqnarray*}
&&(3'\wedge2)_w\\
&=&-m(s'_0,s_0)\cdots
(m-1)(s_{k},s'_{k})(m-1)(s_{k+1},s'_{k+1})(m-1)(t,y_{k+1})\\
&&+(m-1)(s_0,s'_0)(m-1)(s_1,s'_1)\cdots
(m-1)(s_{k+1},s'_{k+1})m(t,y_{k+1})\\
&\equiv&0.
\end{eqnarray*}

Case 2. $w=(m-1)(s_0,s'_0)(m-1)(s_1,s'_1)\cdots
(m-1)(s_k,s'_k)s_{k+1}(2i)(s'_{k+1},s_{k+1})$ $m(s'_{k+1},s_{k+1}),\
0\leq i\leq (m_{s_{k+1}s'_{k+1}}-2)/2$. We consider only the case
that $m_{s_{k+1}s'_{k+1}}$ is odd. The proof is similar for
$m_{s_{k+1}s'_{k+1}}$ to be even. Assume that $m_{s_{k+1}s'_{k+1}}$
is odd. Then
\begin{eqnarray*}
&&(3'\wedge2)_w\\
&=&-m(s'_0,s_0)\cdots
(m-1)(s_{k},s'_{k})(m-1)(s_{k+1},s'_{k+1})(1+2i)(s'_{k+1},s_{k+1})\\
&&+(m-1)(s_0,s'_0)(m-1)(s_1,s'_1)\cdots
(m-1)(s_k,s'_k)s_{k+1}(2i)(s'_{k+1},s_{k+1})m(s_{k+1},s'_{k+1})\\
&\equiv&-m(s'_0,s_0)(m-1)(s_1,s'_1)\cdots
(m-1)(s_{k-1},s'_{k})(m-2-2i)(s_{k+1},s'_{k+1})\\
&&+(m-1)(s_0,s'_0)(m-1)(s_1,s'_1)\cdots(m-1)(s_{k},s'_{k})(m-1-2i)(s'_{k+1},s_{k+1})\\
&\equiv&-(m-1)(s_0,s'_0)(m-1)(s_1,s'_1)\cdots
m(s_{k},s'_{k})(m-2-2i)(s_{k+1},s'_{k+1})\\
&&+(m-1)(s_0,s'_0)(m-1)(s_1,s'_1)\cdots(m-1)(s_{k},s'_{k})(m-1-2i)(s'_{k+1},s_{k+1})\\
&\equiv&-(m-1)(s_0,s_0)(m-1)(s_1,s'_1)\cdots
(m-1)(s_{k},s'_{k})(m-1-2i)(s'_{k+1},s_{k+1})\\
&&+(m-1)(s_0,s'_0)(m-1)(s_1,s'_1)\cdots(m-1)(s_{k},s'_{k})(m-1-2i)(s'_{k+1},s_{k+1})\\
&\equiv&0.
\end{eqnarray*}

\item[($3'\wedge3'$)]\ There are two cases to consider.

Case 1. $w=(m-1)(s_0,s'_0)(m-1)(s_1,s'_1)\cdots
(m-1)(s_k,s'_k)m(s_{k+1},s'_{k+1})
(m-2)(t,y_{k+1})(m-1)(t_1,t'_1)\cdots
(m-1)(t_l,t'_l)m(t'_{l+1},t_{l+1}), \ \ y_{k+1}>t,$ where $y_{k+1}$
is the last letter of $m(s_{k+1},s'_{k+1})$.
\begin{eqnarray*}
&&(3'\wedge3')_w\\
&=&-m(s'_0,s_0)\cdots (m-1)(s_{k+1},s'_{k+1})\cdot(m-2)(t,y_{k+1})
(m-1)(t_1,t'_1)\cdots
m(t'_{l+1},t_{l+1})\\
&&+(m-1)(s_0,s'_0)\cdots (m-1)(s_{k+1},s'_{k+1})m(t,y_{k+1})\cdots
(m-1)(t_{l+1},t'_{l+1})\\
&\equiv&-m(s'_0,s_0)\cdots(m-1)(s_{k+1},s'_{k+1})(m-1)(t,y_{k+1})\cdots
(m-1)(t_{l+1},t'_{l+1})\\
&&+m(s'_0,s_0)\cdots(m-1)(s_{k+1},s'_{k+1})(m-1)(t,y_{k+1})\cdots(m-1)(t_{l+1},t'_{l+1})\\
&\equiv&0.
\end{eqnarray*}

Case 2. $w=(m-1)(s_0,s'_0)(m-1)(s_1,s'_1)\cdots
(m-1)(s_k,s'_k)s_{k+1}(2i)(s'_{k+1},s_{k+1})
(m-1)(s'_{k+1},s_{k+1})(m-1)(t_1,t'_1)\cdots
(m-1)(t_l,t'_l)m(t_{l+1},t'_{l+1}), \ \ \ 0\leq i\leq
(m_{s_{k+1}s'_{k+1}}-2)/2$. We only consider the case that
$m_{s_{k+1},s'_{k+1}}$ is odd and the proof is similar for
$m_{s_{k+1},s'_{k+1}}$ to be even. Assume that
$m_{s_{k+1},s'_{k+1}}$ is odd. Then
\begin{eqnarray*}
&&(3'\wedge3')_w\\
&=&-m(s'_0,s_0)(m-1)(s_1,s'_1)\cdots
(m-1)(s_{k},s'_{k})(m-1)(s_{k+1},s'_{k+1})\\
&&\ \ \cdot (2i)(s'_{k+1},s_{k+1})(m-1)(t_1,t'_1)\cdots
(m-1)(t_l,t'_l)m(t_{l+1},t'_{l+1})\\
&&+(m-1)(s_0,s'_0)(m-1)(s_1,s'_1)\cdots
(m-1)(s_{k},s'_{k})s_{k+1}(2i)(s'_{k+1},s_{k+1})\\
&&\ \ \cdot m(s_{k+1},s'_{k+1}) (m-1)(t_1,t'_1)\cdots
(m-1)(t_{l+1},t'_{l+1})\\
&\equiv&-m(s'_0,s_0)(m-1)(s_1,s'_1)\cdots
(m-1)(s_{k},s'_{k})(m-1-2i)(s_{k+1},s'_{k+1})\\
&&\ \ \cdot (m-1)(t_1,t'_1)\cdots(m-1)(t_{l},t'_{l})m(t_{l+1},t'_{l+1})\\
&&+(m-1)(s_0,s'_0)(m-1)(s_1,s'_1)\cdots
(m-1)(s_{k},s'_{k})(m-1-2i)(s'_{k+1},s_{k+1})\\
&&\ \ \cdot (m-1)(t_1,t'_1)\cdots
(m-1)(t_{l+1},t'_{l+1})\\
&\equiv&-m(s'_0,s_0)(m-1)(s_1,s'_1)\cdots
(m-1)(s_{k},s'_{k})(m-2-2i)(s_{k+1},s'_{k+1})\\
&&\ \ \cdot (m-1)(t_1,t'_1)\cdots(m-1)(t_{l},t'_{l})(m-1)(t_{l+1},t'_{l+1})\\
&&+(m-1)(s_0,s'_0)(m-1)(s_1,s'_1)\cdots
(m-1)(s_{k},s'_{k})(m-1-2i)(s'_{k+1},s_{k+1})\\
&& \ \ \cdot (m-1)(t_1,t'_1)\cdots
(m-1)(t_{l+1},t'_{l+1})\\
&\equiv&-(m-1)(s_0,s'_0)(m-1)(s_1,s'_1)\cdots
(m-1)(s_{k},s'_{k})(m-1-2i)(s'_{k+1},s_{k+1})\\
&&\ \ \cdot (m-1)(t_1,t'_1)\cdots(m-1)(t_{l},t'_{l})(m-1)(t_{l+1},t'_{l+1})\\
&&+(m-1)(s_0,s'_0)(m-1)(s_1,s'_1)\cdots
(m-1)(s_{k},s'_{k})(m-1-2i)(s'_{k+1},s_{k+1})\\
&& \ \ \cdot (m-1)(t_1,t'_1)\cdots
(m-1)(t_{l+1},t'_{l+1})\\
&\equiv&0.
\end{eqnarray*}

\end{enumerate}

Thus, the theorem is proved. \hfill $\blacksquare$

 \ \

We give some examples which are in the case of Theorem \ref{t3.19}
but not the finite Coxeter groups (see \cite{bs01, Lee, Sv}).

\begin{example}
Let $W$ be the Coxeter group with respect to Coxeter matrix
$M=(m_{ij})$. Suppose that one of the following conditions holds:
\begin{enumerate}
\item[(i)]\ for any $i,j\ (i>j)$, $ m_{ij}\geq 3$;

\item[(ii)]\ for any $i,j\ (i>j)$,
either $ m_{ij}=2$ or $m_{ij}=\infty$;

\item[(iii)]\  $ m_{i1}=2$ for any $i\geq 2$ and $m_{ij}\geq 3$
for any $i,j\ (i>j\geq2)$.
\end{enumerate}
Then in $(3')$, there are no relations with property ${\bf C1}\vee
{\bf C2}\vee {\bf C3}\vee {\bf C4}$. By Theorem \ref{t3.19},
$S=\{(\ref{e1}),(\ref{e2}),(3')\}$ is a Gr\"{o}bner-Shirshov basis
of such a Coxeter group $W$.
\end{example}

 \ \

In the next paper, we will try to prove that the new conjecture is
true if $W$ is a Coxeter group without ${\bf C2}\vee {\bf C3}\vee
{\bf C4}$.

\ \

\noindent{\bf Acknowledgement}: The authors would like to thank
Professor L.A. Bokut for his guidance, useful discussions and
enthusiastic encouragement in writing up this paper.

\end{document}